\documentclass[preprint,authoryear]{elsarticle}
\usepackage{graphicx}
\usepackage{lscape}
\usepackage{amsthm}
\usepackage{amsfonts}
\usepackage{amssymb}
\usepackage{dsfont}
\usepackage{amsmath,times,bm,bbm}
\usepackage{xcolor,ulem}

\newtheorem{theorem}{Theorem}

\newdefinition{rmk}{Remark}
\newproof{pf}{Proof}
\newproof{pot}{Proof of Theorem \ref{thm2}}
\newtheorem{corol}{Corollary}

\begin{document}
\bibliographystyle{elsarticle-harv}

\title{Multi-sample rank tests for location against Lehmann-type alternatives}

\author[jeni,nick]{Nikolay I. Nikolov}
\ead{n.nikolov@math.bas.bg}

\author[jeni]{Eugenia Stoimenova\corref{cor1}}
\ead{jeni@math.bas.bg}


\cortext[cor1]{Corresponding author.}

\address[jeni]{Institute of Mathematics and Informatics, Bulgarian Academy of Sciences,
Acad. G.Bontchev str., block 8, 1113 Sofia, Bulgaria}

\address[nick]{Sofia University “St. Kliment Ohridski”, Sofia 1504, Bulgaria}


\begin{abstract}
This paper deals with testing  the equality of $k$  ($k\ge 2$) distribution functions against possible stochastic ordering among them. Two classes of rank tests are proposed  for this testing problem. The statistics of the tests under study are based on precedence and exceedance statistics and are natural extension of corresponding statistics for the two-sample testing problem.  Furthermore, as an extension of the Lehmann alternative for the two-sample location problem, we propose a new subclass of the general alternative for the stochastic order of multiple samples. We show that under the new Lehmann-type alternative any rank test statistics is distribution free. The power functions of the two new families of  rank tests are compared to the power performance of the Jonckheere-Terpstra rank test.

\end{abstract}


\maketitle

\section{Introduction. $k$-sample testing problem}

\noindent{\bf $k$-sample hypothesis testing problem}. Let $X_{i,1}, \ldots, X_{i,n_i}$, \ $i = 1,\ldots ,k$, be $k$ independent random samples with $X_{i,j}$, $j = 1, \ldots , n_i$, having absolutely continuous distribution function $F_i (x)$. We wish to test the hypothesis $H_0$ of homogeneity  against the ordered alternative:
\begin{eqnarray}
H_0: && F_1(x)\equiv \ldots \equiv F_k(x)  \nonumber\\
H_A: &&F_1(x)\ge \ldots \ge F_k(x), \label{genalt}
\end{eqnarray}
where each inequality is strict for some $x$. Note that the distribution of any rank statistic under $H_0$ is not affected by the underlying distribution function, here $ F_1(x)$. However, due to the great generality of the alternative distributions in $H_A$, the distribution under the alternative and some comparative study of the power of tests  could be studied only for particular subclasses of $H_A$.

Ordered alternatives appear to be natural in experiments involving increased levels of a treatment. The most common approach to expressing stochastic dominance in the alternative hypothesis is the location shift parameters alternative, which is in the form $H_{\theta}: F_i(x) =  F(x-\theta_i)$, with $\theta_1\le \ldots \le \theta_k$. However, under such alternatives the distribution of a rank statistic depends not only on $\theta$, but also on $F$. Moreover, parametric tests perform poorly when the underlying distributions are contaminated by possible outliers. Thus, nonparametric settings on the alternative and the tests are particularly important.

In order to study the properties of a nonparametric procedure, it is essential the alternatives themselves to be given in a nonparametric form. The problem of testing homogeneity against ordered alternatives in the multi-sample setup was first considered  by \cite{Terpstra:52} and  \cite{Jonckheere:54}. They suggested the nonparametric test (JT) based on a sum of the Mann-Whitney $U$-statistics associated with the $k(k-1)/2$ possible pairs between the $k$ ordinal groups.

In this paper, we propose a subclass of the alternative $H_A$ which is an extension of the alternatives proposed originally by \cite{lehmann53} for the two-sample testing problem. Lehmann proposed the alternative distributions $H_{L_{2,\eta}}: F_2(x)=\left[F_1(x)\right]^{\eta}$  for some  $\eta \ge 1$. Under $H_{L_{2,\eta}}$ the distribution of any rank statistic does not depend on the underlying distribution functions $F_1(x)$ and $F_2(x)$. Note that, if $\eta$ is a positive integer in this alternative, the second distribution, $F_2(x)$, is indeed the distribution of the largest order statistic from a sample with size $\eta$ from the distribution $F_1(x)$. We extend the idea of using order statistics in the alternatives for comparison of more than two distributions.

\bigskip

\noindent{\bf Lehmann-type alternatives for more than two samples}.
Let $Z_{[1]} \le \ldots \le Z_{[k]}$ be the order statistics of a set of independent observations  $Z_{1},\ldots,Z_{k}$, each having the same distribution function $F(x)$. Define

\begin{equation}\label{LE1}
 H_{L}: \left\{
\begin{array}{lcl}
F_1^*(x)  &=& P(Z_{[1]} \le x)= 1-(1-F(x))^{k}\\
\cdots && \cdots \\
F_i^*(x)  &=& P(Z_{[i]} \le x)= \sum_{t=i}^k \binom{k}{t} [F(x)]^t [1-F(x)]^{k-t}\\
\cdots && \cdots \\
F_k^*(x)  &=& P(Z_{[k]} \le x)= F(x)^{k}.
\end{array}
\right.
\end{equation}

Clearly $F_1^*(x)\ge \ldots \ge F_k^*(x)$ and hence $H_{L}$ is a subclass of $H_A$, defined in (\ref{genalt}). Moreover, it can be used instead of location shift parameters alternative $H_{\theta}$.

The aim of this paper is to investigate rank tests based on precedence and exceedance statistics  for the $k$-sample homogeneity testing problem. We propose two families of statistics that are built as extensions of the corresponding test procedures for the two-sample problem.  Furthermore, we show that the suggested rank tests   are distribution free under the alternative $H_{L}$  in \eqref{LE1}.

\bigskip

The very early results for Lehmann alternatives are due to  \cite{lehmann53} and \cite{savage1956contributions}. These and many other consequent papers are dealing with testing  the equality of two distribution functions against stochastic ordering among them. Some basic early references include \cite{shorack1967tables}, \cite{Davies:71},  \cite{peto1972rank},  \cite{lin1990powers}, \cite{sukhatme1992powers}. A collection of articles devoted to precedence-type tests for the two-sample testing problem  under Lehmann alternative includes \cite{katzen:85}, \cite{vanderLaan:01, van1999best} among  others. One  key feature of all these articles is that the tests are based only on the ranks which makes them permutation invariant in some sense. Usually the main focus is in comparison between the power functions of rank tests under the Lehmann alternatives. \cite{balakrishnan2006precedence} proposed a number of rank tests based on precedence statistics, which they have studied under the null hypothesis and Lehmann alternative. Regarding the multi-sample nonparametric tests, we should mention the papers by \cite{govindarajulu1971c}, \cite{kossler2000efficacy},  \cite{kossler2005some}, \cite{altunkaynak2021comparing}, and the work by \cite{vock2011jonckheere} that considers multi-sample test procedures in the framework of balanced ranked set sampling schemes.

\bigskip

The rest of this paper is organized as follows. In section~\ref{Xdistr}, we derive the distributions of the ordered observations under the new Lehmann-type alternatives. In particular, we describe the connection with distribution of ranks in section \ref{Rconnect}, whereas the probability results for the ordered observations are given in sections \ref{rankDistribution} and \ref{6typesAlternatives}. In section \ref{3tests}, we describe two classes of rank test for the equality of $k$ distribution functions against possible stochastic ordering among them. The rank tests under study are based on precedence and exceedance statistics. In section \ref{3thenull}, tables with critical values for different number of samples and for some selected choices of sample sizes are presented for the new proposed tests. Finally, in section \ref{power} we compare the power functions of the two new families to the Jonckheere-Terpstra rank test, while in section~\ref{additionalRemarks} we provide some additional concluding remarks.

\section{Distributions under new Lehmann-type alternatives}\label{Xdistr}

Before deriving the distribution under the Lehmann type alternatives in \eqref{LE1}, we briefly present the relation between the order statistics in $k$ samples and their corresponding ranks by following the exposition in \cite{hoeffding1951optimum}.

\subsection{Connection with distribution of ranks}\label{Rconnect}

Let ${\bf x}= \left(x_{1,1}, \ldots, x_{1,n_1},\  x_{2,1}, \ldots, x_{2,n_2},  \ldots, \ x_{k,1}, \ldots, x_{k,n_k}\right)\in \mathcal{X}$ be the obtained observations from $k$ independent samples, where $\mathcal{X}$ is a subset in a Euclidean space of dimension $n=\sum_{i=1}^{k}n_i$. Furthermore, we assume that ${\bf x}$ is a realization of a random vector ${\bf X}=( X_{1,1}, \ldots, X_{1,n_1}$,\  $X_{2,1}, \ldots, X_{2,n_2}$,  $\ldots, X_{k,1}, \ldots, X_{k,n_k})$ with  continuous distribution. Thus, with probability 1 there are no ties in ${\bf x}$, i.e. ${\bf x}\in W$, where $W$ is the subset of $\mathcal{X}$ such that $x_{i,j}\neq x_{g,h}$ if $i\neq g$ or $j\neq h$. We say that the coordinate $x_{i,j}$ has rank $r_{i,j}$ with respect to ${\bf x}$ if exactly $r_{i,j}-1$ coordinates $x_{g,h}$ are less than $x_{i,j}$, for $g\neq i$ and $h\neq j$. Let us denote by ${\bf r}$ the rank vector associated with ${\bf x}$, i.e.
\begin{equation*}
	{\bf r} =(r_{1,1}, \ldots, r_{1,n_1},\  r_{2,1}, \ldots, r_{2,n_2},   \ldots, \ r_{k,1}, \ldots, r_{k,n_k}).
\end{equation*}
Since ${\bf x}\in W$ has no ties, ${\bf r}$ is uniquely determined by ${\bf x}$, with ${\bf r}$ being a permutation of $\left\{1,2\ldots,n\right\}$. Furthermore, a fixed rank permutation ${\bf r}$ is associated with a subset of $W$. Thus, the set $W$ can be expressed as a union of $n!$ non-intersecting subsets, each corresponding to a rank permutation. Therefore, in the nonparametric testing approach based on rank order the critical region for rejecting the null hypothesis $H_{0}$ is a union of such subsets that are associated with the alternative hypothesis (rank permutations under the alternative). Moreover, since the observations in each of the $k$ samples are independent and identically distributed, the observed ranks are invariant in permutations of the observations within each sample. Hence, instead of considering the $n!$ possible rankings of all observations, we can study only $\frac{n!}{n_1!\ldots n_k!}$ groups of rankings with $n_1!\ldots n_k!$ permutations in every group. In such way the test procedures are based on the order statistics of the $k$ samples. A rigorous formulation of rank test methods is given in the works of \cite{lehmann53} and  \cite{hoeffding1951optimum}, while for more properties of the nonparametric rank tests we refer to the books of \cite{lehmann2022testing} and \cite{sidak1999theory}.

\subsection{Distribution of ordered observations under the Lehmann type alternatives}\label{rankDistribution}
In order to obtain the distribution of the ordered observations under the Lehmann type alternatives we need the following notations. As in Section~\ref{Rconnect}, let ${\bf x}= (x_{1,1}, \ldots, x_{1,n_1},$ $ \ldots, \ x_{k,1}, \ldots, x_{k,n_k})$ and ${\bf r}= \left(r_{1,1}, \ldots, r_{1,n_1}, \ldots, \ r_{k,1}, \ldots, r_{k,n_k}\right)$ be the vector of observations and associated ranks, respectively, for $k$ independent samples from continuous random variables. Denote by $a_{j}({\bf r})$ the sample index of the observation in ${\bf x}$ that has rank $j$ in ${\bf r}$, for $j=1,2,\ldots,n$, i.e.,
\begin{equation*}\label{defA}
	a_{j}({\bf r})=\sum_{i=1}^{k}\sum_{h=1}^{n_{i}}i\mathbbm{1}\left\{r_{i,h}=j\right\},
\end{equation*}
where $\mathbbm{1}\left\{\cdot\right\}$ is the indicator function. Furthermore, for the observation with rank $j$ let $b_{j}({\bf r})$ be the index within its sample, i.e.,
\begin{equation*}\label{defB}
	b_{j}({\bf r})=\sum_{h=1}^{n_{a_{j}({\bf r})}}h\mathbbm{1}\left\{r_{a_{j}({\bf r}),h}=j\right\}, \quad \mbox{for $j=1,2,\ldots,n$.}
\end{equation*}
The following theorem provides an expression for the probability of observing a fixed rank order under the Lehmann-type alternative $H_{L}$ in \eqref{LE1}.
\begin{theorem}\label{ThOrder}
	Let ${\bf X}= (X_{1,1}, \ldots, X_{1,n_1}$,\  $X_{2,1}, \ldots, X_{2,n_2}$,  $\ldots, X_{k,1}, \ldots, X_{k,n_k})$ be $k$ samples of mutually independent random variables such that $X_{i,j}$ ($i=1, \ldots, k$, $j=1, \ldots, n_i$) has distribution function
	\begin{equation}\label{LE}
		H_{L}: F_{i,k}(x)  =  \sum_{t=i}^k \binom{k}{t} [F(x)]^t [1-F(x)]^{k-t}.
	\end{equation}
	where $F(x)$ is an arbitrary unknown continuous distribution function. Then, for the rank order vector ${\bf R}= \left(R_{1,1}, \ldots, R_{1,n_1}, \ldots, \ R_{k,1}, \ldots, R_{k,n_k}\right)$ of ${\bf X}$ and given permutation ${\bf r}= (r_{1,1}, \ldots, r_{1,n_1}, \ldots,$ $r_{k,1}, \ldots, r_{k,n_k})$ of $\left\{1,2\ldots,n\right\}$,
	\begin{align}
		\nonumber &P({\bf R}={\bf r}| H_{L}) = P(X_{a_{1}({\bf r}),b_{1}({\bf r})}\leq X_{a_{2}({\bf r}),b_{2}({\bf r})}\leq \ldots \leq X_{a_{n}({\bf r}),b_{n}({\bf r})}| H_{L})\\ &=
		 \idotsint\limits_{0 \le y_{a_{1}({\bf r}),b_{1}({\bf r})}\leq y_{a_{2}({\bf r}),b_{2}({\bf r})}\leq \ldots \leq y_{a_{n}({\bf r}),b_{n}({\bf r})} \le 1} \ \prod_{i=1}^{k}\prod_{j=1}^{n_i} d\left( \sum_{t=i}^k \binom{k}{t} y_{i,j}^t (1-y_{i,j})^{k-t}  \right), \label{ProbOrder}
	\end{align}
	
\end{theorem}
\begin{proof}
	Since the variables in ${\bf X}$ are mutually independent with distribution functions in \eqref{LE}, it follows that
	\begin{align}\label{integralLE}
		\nonumber &P(X_{a_{1}({\bf r}),b_{1}({\bf r})}\leq X_{a_{2}({\bf r}),b_{2}({\bf r})}\leq \ldots \leq X_{a_{n}({\bf r}),b_{n}({\bf r})}| H_{L})=\\
		&\idotsint\limits_{-\infty < x_{a_{1}({\bf r}),b_{1}({\bf r})}\leq x_{a_{2}({\bf r}),b_{2}({\bf r})}\leq \ldots \leq x_{a_{n}({\bf r}),b_{n}({\bf r})} < \infty} \ \prod_{i=1}^{k}\prod_{j=1}^{n_i} d\left( \sum_{t=i}^k \binom{k}{t} [F(x_{i,j})]^t [1-F(x_{i,j})]^{k-t}  \right).
	\end{align}
	The fact that $F(x)$ is a continuous distribution function implies that $0\leq F(x) \leq 1$ and its inverse $F^{-1}(x)$ is a non-decreasing continuous function. Thus, \eqref{ProbOrder} is obtained by changing the variables in \eqref{integralLE} as $x_{i,j}=F^{-1}(y_{i,j})$, for $i=1,\ldots,k$ and $j=1,\ldots,n_i$.
\end{proof}

Note that the probability in \eqref{ProbOrder} does not depend on the underlying unknown distribution $F(\cdot)$, but only on the ordering in ${\bf X}$ through the rank permutation ${\bf r}$ and the indices  $a_{j}({\bf r})$ and $b_{j}({\bf r})$, for $j=1,2\ldots,n$. Furthermore, as the observations are invariant in permutations within each sample, the probabilities for rankings of the order statistics are obtained by multiplying the probability \eqref{ProbOrder}, for the corresponding rank ${\bf r}$, by $n_1!\ldots n_k!$.


If we consider the special case of the most extremal ordering, associated with the alternative $H_{A}$ in \eqref{genalt}, the following result is directly implied by Theorem~\ref{ThOrder}.
	\begin{corol}\label{corollary3}
		Under the same notations as in Theorem~\ref{ThOrder}, for the identity permutation $(1,\ldots, n)$, we obtain that
		\begin{align*}
			\nonumber P(X_{1,1}^{(1)} &\le \ldots \le X_{1,n_1}^{(1)} \le \ldots \le X_{k,1}^{(k)} \le \ldots \le X_{k,n_k}^{(k)}| H_{L}) = \\
			&  \idotsint\limits_{0 \le  y_{1,1}\leq y_{1,2}\leq \ldots \leq y_{k,n_{k}} \le  1} \ \prod_{i=1}^{k}\prod_{j=1}^{n_i} d\left( \sum_{t=i}^k \binom{k}{t} y_{i,j}^t (1-y_{i,j})^{k-t}  \right).
		\end{align*}
	\end{corol}


\subsection{Other types of Lehmann's alternatives}\label{6typesAlternatives}

The original Lehmann alternative introduced by  \cite{lehmann53} is for testing equality of the distributions with hypothesis $H_0: F(x)= G(x)$ against $H_{L_{2,\eta}}:G(x) = F^{\eta}(x)$ ($\eta>1$). Obviously, $H_{L_{2,\eta}}$ is a subclass of the general alternative $H_{A}$. Lehmann showed that under $H_{L_{2,\eta}}$ the distribution of any rank statistic is distribution free. \cite{savage1956contributions} considers a  more generalized form of the Lehmann alternative and shows the following results, presented in Theorem~\ref{SavageTh} and Corollary~\ref{SavageCor} below.

\begin{theorem}[Savage, 1956]\label{SavageTh}
	
	Let the random variables $X_1, \ldots, X_N$ be mutually independent such that $X_i$ has the distribution function $H_{L_{\eta}}: F_i(x) = [H(x)]^{\eta_i}$ , $i = 1, \ldots N$, where $\eta_i > 0$, and $H(x)$ is an unknown continuous distribution function. Then
	\begin{equation}\label{SavageFormula}
	P(X_{k_1} \le X_{k_2} \le  \cdots \le X_{k_N} \ | H_{L_{\eta}}) = \left( \prod_{i=1}^{N} \eta_i \right) / \prod_{i=1}^{N} \left( \sum_{j=1}^{i}   \eta_{k_j} \right),
	\end{equation}
	where $(k_1,\ldots, k_N)$ is a permutation of
	$\left\{1,2\ldots,N\right\}$
	and $H_{L_{\eta}}$ denotes the specified Lehmann alternative.
\end{theorem}

\begin{corol}\label{SavageCor}
	Under the same notations as in Theorem~\ref{SavageTh}, for the identity permutation $(1,\ldots, N)$ it holds that
	\begin{equation*}\label{SavageCorFormula}
	P(X_{1} \le X_{2} \le  \cdots \le X_{N} \ | H_{L_{\eta}}) = \left( \prod_{i=1}^{N}\eta_i \right) / \prod_{i=1}^{N} \left( \sum_{j=1}^{i}  \eta_{j} \right).
	\end{equation*}

\end{corol}

Here, we modify Salvage's theorem for the $k$-sample problem for testing the hypothesis in \eqref{genalt}.

	\begin{theorem}\label{SavageKsampleTh}
		Let the random variable $X_{i,j}$, $j= 1, \ldots n_i$, from the $i$-th sample, $i = 1,\ldots ,k$, has distribution function
		\begin{equation}\label{Leta}
			H_{L_{k,\eta}}: F_i(x) = [F(x)]^{\eta_i},
		\end{equation}
		 where $\eta_1 \le \ldots \le \eta_k$. Then, for  a given permutation ${\bf r}= (r_{1,1}, \ldots, r_{1,n_1}, \ldots,$ $r_{k,1}, \ldots, r_{k,n_k})$ of $\left\{1,2\ldots,n\right\}$,
		
		\begin{align}
			\nonumber P({\bf R}={\bf r}|  H_{L_{k,\eta}}) = &P(X_{a_{1}({\bf r}),b_{1}({\bf r})}^{(a_{1}({\bf r}))}\leq X_{a_{2}({\bf r}),b_{2}({\bf r})}^{(a_{2}({\bf r}))}\leq \ldots \leq X_{a_{n}({\bf r}),b_{n}({\bf r})}^{(a_{n}({\bf r}))}| H_{L_{k,\eta}})\\ =&
			\left( \prod_{i=1}^{k} n_{i}!\eta_i^{n_i} \right) / \prod_{i=1}^{n} \left(\sum_{j=1}^{i} \eta_{a_{j}({\bf r})}\right), \label{SavageGeneralFormula}
		\end{align}
		where $X_{i,1}^{(i)} \le \ldots \le X_{i,n_i}^{(i)}$, $i=1, \ldots, k$, denote the ordered observations in each sample.
		
	\end{theorem}	

	\begin{proof}
		 Formula~\eqref{SavageGeneralFormula} follows directly from \eqref{SavageFormula} in Theorem~\ref{SavageTh}	and the fact that the observations are invariant in permutations within each sample.
	\end{proof}
If we consider the special case of the most extremal ordering, associated with the alternative $H_{A}$ in \eqref{genalt}, we obtain the following result as a straightforward implication of Theorem~\ref{SavageKsampleTh} for the identity permutation $(1,\ldots, n)$.
	\begin{corol}\label{corollary2}
		  Let the same notations as in Theorem~\ref{SavageKsampleTh} hold. Then,
		\begin{align*}
		P(X_{1,1}^{(1)} \le \cdots \le X_{1,n_1}^{(1)}\le &\cdots \le X_{k,1}^{(k)} \le \cdots \le X_{k,n_k}^{(k)} \ | H_{L_{k,\eta}}) \\& = \left( \prod_{i=1}^{k} n_{i}! \eta_i^{n_i} \right) / \prod_{i=1}^{k} \prod_{j=1}^{n_{i}}
		\left( j\eta_{i}+\sum_{l=1}^{i-1} n_{l}\eta_{l}   \right).
		\end{align*}
	\end{corol}

	The probabilities for the identity rank permutation in Corollaries~\ref{corollary3}, \ref{SavageCor} and \ref{corollary2} play an important role for rejecting the null hypothesis in the nonparametric rank approach for testing \eqref{genalt}. In the next section, we present an unified method for developing rank tests by constructing an extremal set of ranks, which are similar to the identity rank.

\section{Rank data structure for multi-sample location problem}\label{structure}

We aim to develop a nonparametric test for the problem \eqref{genalt}, on the basis of  the observed rank vector. The approach is based on the distance between two sets of permutations and was first proposed by \cite{critchlow92}. The general method of constructing a suitable rank test for testing $H_0$ versus $H_A$ is implemented as follows:

Let $X_{i,1}, \ldots, X_{i,n_i}$, \ $i = 1,\ldots ,k$, are $k$ independent random samples with $X_{i,j}$, $j = 1, \ldots , n_i$ having an absolutely continuous distribution function $F_i (x)$. Denote $n=\sum_{j=1}^k n_j$.

{\em Step 1.} Rank all the $n$ observations among themselves and denote the rank of $X_{i,j}$ in the pooled sample by $R_{i,j}$ \  $(i=1,\ldots,k; j=1,\ldots, n_i )$. Thus,
\[
{\bf R} = \left( R_{1,1}, \ldots, R_{1,n_1}, \ldots, R_{k,1}, \ldots, R_{k,n_k}\right),
\]
is a permutation from the set $S_n$ of all permutations of $n$ elements. Note, that the ranks from the same distribution are placed together in ${\bf R}$.

{\em Step 2.} Denote by $[{\bf R}]$  the equivalent subclass consisting all permutations of the integers $1, \ldots , n$ which assign the same set of ranks to the individual populations as ${\bf R}$ does. Thus $[{\bf R}]$ contains $n_1! \ldots n_k!$ permutations, and is the left coset of $S_n$.

{\em Step 3.} Further, identify the set $E$ of extremal permutations, consisting of all permutations that are least in agreement with $H_0$, and most in agreement with $H_A$. The extremal set $E$ consists of all permutations which assign ranks $ c_{j-1} + 1, \ldots, c_j$ to population $j$, where $c_j=\sum_{i=1}^{j} n_i, \ j=1,\ldots, k$.

{\em Step 4.} To test $H_0$ versus $H_A$  let $d$ be a distance on $S_n$ and define the rank test by the distance $d([{\bf R}],E)$.

Since $d([{\bf R}],E)$ measures the distance from $[{\bf R}]$  to the set of permutations that are most in agreement with $H_A$, it follows that $H_0$ should be rejected for small values of $d([{\bf R}],E)$. It is natural to require that the distance between rankings should  be invariant under an arbitrary relabeling in each sample (right-invariant). Distances involved in rank test are usually right-invariant.

One of the key features of the tests based  on ranks is that they remain invariant under arbitrary relabeling in each sample. Consequently, we may relabel the observations in each sample so that their ranks are arranged in ascending order.
As it was mentioned by \cite{critchlow92}, the method gives rise to most of the standard rank statistics and also leads to some new rank test statistics. Similar approach is used by \cite{alvo1997general} and describes various test statistics for the multi-sample location problem with ordered alternatives as well. However, the distributions of any tests statistics under the hypotheses should be studied in any particular case.

\section{Rank tests for ordered alternatives}\label{3tests}

We will explore the performance of two classes of rank tests for testing the hypotheses  in \eqref{genalt}.

The proposed test statistics  rely on partially observed data sets with the aim to reduce the influence of possible outliers in the underlying population distributions. The tests are extensions of the corresponding tests for the two-sample homogeneity testing, which are studied by \cite{jspi:11, stoimenova2017vsidak}. The new proposed tests will be studied with comparison to the Jonckheere-Terpstra, suggested by \cite{Jonckheere:54} and \cite{Terpstra:52}.

Recall that $X_{i,1}, \ldots, X_{i,n_i}$, \ $i = 1,\ldots ,k$, are $k$ independent random samples with $X_{i,j}$, $j = 1, \ldots , n_i$ having an absolutely continuous distribution function $F_i (x)$. Let $R_{i,j}$ \ be the rank of $X_{i,j}$ in the pooled sample of all observations $(i=1,\ldots,k; j=1,\ldots, n_i)$. Furthermore, let us consider the following key ranks for each sample.

\noindent{\bf Notation.} Denote $R^*_{1,1}\le \cdots \le R^*_{1,n_1}$  the ranks of the first sample (in the pooled sample) arranged in ascending order, $\ldots$, $R^*_{k,1} \le \cdots \le R^*_{k,n_k}$  the ranks of the last sample (in the pooled sample) arranged in ascending order.

\subsection{A class of tests based on maximum deviations}\label{Mtets}

For two samples $X_1,\ldots,X_{m}$  $\sim$ $F$ and $Y_1,\ldots,Y_{n}$  $\sim$ $G$ the $M-$statistic is defined by $M  = \max\{R^*_{1,m} - m, m+1 - R^*_{2,1}\}$, where $R^*_{1,m}$ is the largest rank assigned to the $X$-sample and $R^*_{2,1}$ is the smallest rank assigned to the $Y$-sample in the pooled sample. The test is induced by the Critchlow's method, described in section~\ref{structure}, with the Chebyshev metric on the set of all permutations.

The statistic $M$ is appropriate for testing the hypothesis that $F$ and $G$ are identical against the one-sided alternative that $Y$'s are stochastically larger than $X$'s.
\cite{comm:11} studied in more detail the case with equal sample sizes.  Similar tests, based on extreme ranks of the observations, have been proposed by \cite{haga1959two}, \cite{sidak:57}, \cite{sidak1999theory}, and others.
 \cite{jspi:11} has generalized the $M$-statistic  for two samples as
\begin{equation}\label{Mrho2}
M_{\rho}  = \max\{R^*_{1,m-s} - m+s, m+r+1 -R^*_{2,1+r}\}
\end{equation}
where $s=[\rho m]$ and $r=[\rho n]$ for some {\small $0\le \rho <1$}, with $[ \cdot ]$ denoting  the integer part. The  motivation of introducing the corrections by $s$ and $r$ is to reduce the influence of the extremal ranks $R^*_{1,m}$ and $R^*_{2,1}$ in the underlying observed samples.

Now, we define an extension of the $M_{\rho}$ statistics for the $k$-sample problem, given in \eqref{genalt}.

Under the notations in this section, the extension of the test statistic in \eqref{Mrho2} for the multi-sample problem is given by
\begin{equation}\label{Mrhom}
M_{\rho}  = \max_{1\le j \le k-1}\left[\max \left\{ |R^*_{j,n_j-s_j}-c_{j}+s_j|, \ |c_j+1+s_{j+1}-R^*_{j+1,1+s_{j+1}}|
\right\}\right],
\end{equation}
 where  $s_j=[\rho n_j]$   for some $0\le \rho <1$ and $c_j=\sum_{i=1}^{j} n_i, \ j=1,\ldots, k-1$.

Evidently, small values of $M_{\rho}$ lead to the rejection of $H_0$ in favor of the stochastically ordered alternative $H_A$ in \eqref{genalt}. The test statistic $M_{0}$, obtained for $\rho=0$ in \eqref{Mrhom}, is based on the maximum deviation of the largest and the smallest rank in each sample to the  corresponding ranks in the extremal ranking  $E$, defined in Step~3 in section~\ref{structure}. Various values of $\rho$ yield a family of test statistics which we refer to as $M$-type tests. For $\rho>0$ we expect the tests to perform better against outliers, as extreme ranks from the samples are not used when the corrections $s_j$, for $j=1,\ldots,k-1$, are strictly positive. It is reasonable  the value of the trimmed proportion $\rho$ to be non-zero,  since we   desire  to reduce the possible influence of a small number of potential outliers.


\subsection{\v{S}id\'ak-type tests}

For testing the equality of two distribution functions, $F=G$,  \cite{stoimenova2017vsidak} proposed  the following class of rank tests based on precedence and exceedance statistics. Given two samples $X_1,\ldots,X_{m}$  $\sim$ $F$ and $Y_1,\ldots,Y_{n}$  $\sim$ $G$, and under the notations in section~\ref{Mtets}, the suggested test statistic $V_{\rho}$ has the form
\begin{equation}\label{vrstat}
V_{\rho}  =  (R^*_{1,m-s} - m+s) + ( m+r+1 -R^*_{2,1+r}),
\end{equation}
where  $R^*_{1,m-s}$ is the $(m-s)$th largest rank assigned to the $X$-sample, while $R^*_{2,r+1}$ is the $(r+1)$th smallest rank assigned to the $Y$-sample.   For $\rho=0$ the test statistic is equivalent to the test statistic proposed by \v{S}id\'ak,  see \cite{sidak:57} and  \cite{sidak1999theory}. The test is appropriate for testing  the hypothesis that $F$ and $G$ are identical against the one-sided alternative that $Y$'s are stochastically larger than $X$'s.


Now, we define an extension of the $V_{\rho}$ statistics for the multi-sample testing problem   in \eqref{genalt}. Let $X_{i,1}, \ldots, X_{i,n_i} \sim F_j$, $i = 1,\ldots ,k$, be the observed vectors and $R^*_{1,1}\le \cdots \le R^*_{1,n_1}$,    $\ldots$, $R^*_{k,1}\le \cdots \le R^*_{k,n_k}$   be the ordered  ranks in the pooled sample, as in the previous notations in this section.    The first extension we propose is
\begin{equation}\label{V0m}
V_{0} = \sum_{j=1}^{k-1} \left[| R^*_{j,n_j}-c_{j}| +| c_j+1 -R^*_{j+1,1} |\right],
\end{equation}
where $c_j=\sum_{i=1}^{j} n_i, \ j=1,\ldots, k-1$.
Small values of $V_{0}$ lead to the rejection of $H_0$ in favor of the stochastically ordered alternative $H_A$. The test statistic in (\ref{V0m}) is extension of $V_{0}$, defined in (\ref{vrstat}) for $\rho=0$. It is based on the sum of all deviations of the largest and the smallest rank in each sample to the corresponding rank in the extremal ranking set $E$.

Similar to the test statistic $M_\rho$, we make a further generalization of the test statistic  $V_0$ in \eqref{V0m}  and define a class of rank test statistics. The aim is to  reduce the   effect of the  extremal ranks.

Let $s_j=[\rho n_j]$ for some $0\le \rho <1$   and define
\begin{equation}\label{Vrhom}
V_{\rho}  =   \sum_{j=1}^{k-1} \left[ | R^*_{j,n_j-s_j}-c_{j}+s_j| +|  c_j+1+s_{j+1} -R^*_{j+1,1+s_{j+1}}|\right]
\end{equation}
where $c_j=\sum_{i=1}^{j} n_i, \ j=1,\ldots, k-1$.
Small values of $V_{\rho}$ lead to the rejection of $H_0$ in favor of the stochastically ordered alternative $H_A$  in \eqref{genalt}. Various values of $\rho$ yield a family of test statistics, which we refer to as \v{S}id\'ak-type tests.   As the statistic $M_\rho$, the test based on $V_{\rho}$ is expected to be more robust against outliers for strictly positive trimmed proportion $\rho>0$. The role of the parameter $\rho$ on the power of the tests based on $M_\rho$ and $V_\rho$ is investigated via a simulation study in section \ref{power}.

\subsection{Jonckheere-Terpstra test}

The rank-sum statistic  of Jonckheere-Terpstra ($JT$) is known to provide a good nonparametric test for the hypothesis testing problem described above, it will be used for comparison.  Therefore, we will use it for comparison to the two families of statistics $M_\rho$ and $V\rho$.    Under the notations used in this section, let
\begin{equation}\label{JT}
JT  =    \sum_{1 \le i \le j \le k} U_{ij} = \sum_{i=1}^{k-1} \sum_{j=i+1}^{k-1} \sum_{t=1}^{n_i}  \sum_{s=1}^{n_j} I_{\{ X_{i,t} <X_{j,s} \}}
\end{equation}
where $U_{ij}$ is the number of observations in sample $j$ that exceed observations in sample $i$. Thu, $JT$ is a linear sum of several Mann-Whitney statistics, see \cite{Jonckheere:54}.
In terms of ranks, Jonckheere-Terpstra test is equivalent to the sum all deviations from the ordered observed ranks to the identity permutation corresponding to the ideal ranks under the alternative, i.e. to the extremal set $E$. Moreover, the statistic $JT$ can be expressed in the form
\[
JT = \sum_{j=1}^{k} \sum_{t=1}^{n_j}  | R^*_{j,t}-( \sum_{l=1}^{n_j-1}  n_{l}+t   ) |,
\]
with  $H_0$ being  rejected for large values of $JT$.


\section{Null distribution  simulations}\label{3thenull}

Under $H_0$ and for fixed $\rho$, the critical region for $M_{\rho}$ at a pre-specified significance level $\alpha$ is given by $\{ 0, 1,  \ldots, s_{\alpha} \}$, where $s_{\alpha}$ satisfies $ P(M_{\rho}\le s_{\alpha} | H_0)\le \alpha.$
We apply  Monte Carlo simulations to   approximate  $s_{\alpha}$ for particular choice of the number of samples and particular choices of the sizes of the samples. These calculations have been carried out on a PC computer by using the statistical package R.

Since the $M_{\rho}$-statistic is a discrete random variable, the critical value $s_{\alpha}$ is the largest value of $M_{\rho}$ that satisfies $P(M_{\rho}\le s_{\alpha} | H_0)\le \alpha$, where in general we do not obtain equality to $\alpha$. In order to achieve a prescribed level $\alpha$, we may use randomization for each realization of the samples  and  then count the probability of rejecting $H_0$ as follows. Calculate  $s_{\alpha}$ for which
\begin{equation}\label{rand-prp}
  P(M_{\rho}\le s_{\alpha}|H_0) = \alpha_L
   \quad   \mbox{and}   \quad
  P(M_{\rho}\le s_{\alpha}+1|H_0) = \alpha_R,
\end{equation}
where $\alpha_L < \alpha < \alpha_R$. Independent of the given data, produce a Bernoulli random variable $Z$ that takes   value 1 with probability
\begin{equation}\label{pib}
\pi = \frac{\alpha - \alpha_L}{\alpha_R - \alpha_L}
\end{equation}
and value 0 with probability $1-\pi$. Thus, $\pi$  and $1-\pi$ are the proportions of times out of a fixed number of runs that each procedure in (\ref{rand-prp}) rejects the null hypothesis. As a result, the level of significance of a test procedure with critical region
\[
 \{M_{\rho} \le s_{\alpha}\} \cup   (\{M_{\rho} \le s_{\alpha}+1\} \cap \{Z=1\})
\]
will be exactly $\alpha$. Since   this  test procedure achieves the exact level of significance, it would enable us to compare the power of the test for different sample sizes as well as to compare powers of different rank tests.

Under $H_0$, the critical region for $V_{\rho}$ for a pre-specified significance level $\alpha$ has the same form as those of $M_{\rho}$.

In this study the $V_{\rho}$-type tests and $M_{\rho}$-type tests are specified for   $\rho \in \left\{0, 0.05, 0.1, 0.15, 0.20, 0.25\right\}$,  with $s_j$ determined as $s_j = [\rho n_j]$.  Below,  we have presented the critical value (c.v.) at 5\% level of significance, estimated through Monte Carlo simulations based on 10,000 simulated data sets, for each case specified by $\rho$.

Table \ref{cv3table} contains the critical values of the $V_{\rho}$-tests, $M_{\rho}$-tests and Jonckheere-Terpstra test for 3 samples with equal sample sizes $n =5 (5) 25$, and $\rho =0(0.05)0.25$ at 5\% level of significance. The exact  tail probability, without the adjustment by \eqref{rand-prp} and \eqref{pib}, values are given in parentheses (left-tail for the $M_{\rho}$-tests and $V_{\rho}$-tests and right-tail for Jonckheere-Terpstra test).
\begin{table}
	\centering
	\caption{Critical values for 3 samples  with equal sample sizes  $n= 5 (5) 25$, and $\rho =0(0.05)0.25$ at near 5\% level of significance}
	\label{cv3table}
\begin{tabular}{cc||c|c|c}
$\rho$ &$n$   & $V_{\rho}$&    $M_{\rho}$&    $JT$      \\
\hline
     0 &  5    &  18 (0.0345) &  7 (0.0401) &  53 (0.0484) \\
  0.05 &  5    &  18 (0.0345) &  7 (0.0401) &    \\
  0.10 &  5    &  18 (0.0345) & 7  (0.0401)  &    \\
  0.15 &  5    &  18 (0.0345) &  7 (0.0401) &    \\
  0.20 &  5    &  11 (0.0485) &  4 (0.0195)   &    \\
  0.25 &  5    &  11 ( 0.0485) &  4 (0.0195)  &    \\
\hline
     0 &  10   & 46 (0.0347) &  17 (0.0348)   & 192 (0.0479)\\
  0.05 &  10   & 46 (0.0347) &  17 (0.0348)   &    \\
  0.10 &  10   & 38 (0.0490) &  14  (0.0342)  &     \\
  0.15 &  10   & 38 (0.0490) &  14  (0.0342)  &     \\
  0.20 &  10   & 30 (0.0490) & 11  (0.0299)   &    \\
  0.25 &  10   & 30 (0.0490) & 11  (0.0299) &    \\
 \hline
    0 &  15    & 75 (0.0372)   & 54 (0.0402)   &  418 (0.0482)\\
 0.05 &  15    & 75 (0.0372)   & 27 (0.0402) &   \\
 0.10 &  15    &  66 (0.0429)   & 24 (0.0349)  &    \\
 0.15 &  15    &  57 (0.0416)   & 21 (0.0307)  &   \\
 0.20 &  15    &  49 (0.0469)   & 19 (0.0439) &   \\
 0.25 &  15    &  49 (0.0469)   &  19 (0.0439)  &   \\
  \hline
 0 &     20    &  105 (0.0379)  & 37 (0.0339)    & 719 (0.0500) \\
 0.05 &  20    &   95 (0.0370)  & 34 (0.0362 )  &  \\
 0.10 &  20    &   86 (0.0376)  & 31 (0.0305) &    \\
 0.15 &  20    &   78 (0.0451)  & 29 (0.0417) &    \\
 0.20 &  20    &   69 (0.0387)  & 26 (0.0334) &    \\
 0.25 &  20    &   61 (0.0430)  & 24 (0.0389) &    \\
    \hline
   0  &  25    & 135 (0.0447)   & 47 (0.0382) &  1107 (0.0499)\\
 0.05 &  25    & 125 (0.0443)   & 44 (0.0334) &     \\
 0.10 &  25    &  116 (0.0459)   & 41 (0.0303) &     \\
 0.15 &  25    &  107 (0.0441)   &  39 (0.0423) &     \\
 0.20 &  25    & 89 (0.0411)   & 34 (0.0411)  & \\
 0.25 &  25    &  81 (0.0438)  &  32 (0.0468)  &     \\
\hline
\end{tabular}

\end{table}
%
%
%
%
Table \ref{cv4table} contains the critical values of the $V_{\rho}$-tests, $M_{\rho}$-tests and Jonckheere-Terpstra test for 4 samples with equal sample sizes $n=5 (5) 25$, and $\rho =0(0.05)0.25$ at 5\% level of significance. The  exact  tail probabilities are  given in parentheses.
\begin{table}
	\centering
	\caption{Critical values for 4 samples with equal sample sizes $n =5 (5) 25$, and $\rho =0(0.05) 0.25$ at near 5\% level of significance}
	\label{cv4table}
\begin{tabular}{cc||c|c|c}
$\rho$ &$n$   & $V_{\rho}$&    $M_{\rho}$&    $JT$      \\
\hline
     0 &  5    &  38 (0.0442) &  11 (0.047) & 100 (0.0461) \\
  0.05 &  5    &  38 (0.0442) &   11 (0.048) &    \\
  0.10 &  5    &  38 (0.0442) &   11 (0.048) & \\
  0.15 &  5    &  38 (0.0442) &   11 (0.048) &    \\
  0.20 &  5    &  24 (0.0447) &  8 (0.0468)   &    \\
  0.25 &  5    &   24 (0.0447) &  8 (0.0468)   & \\
\hline
     0 &  10   &  94 (0.0439) &  25 (0.028)   & 369 (0.047)\\
  0.05 &  10   &  94 (0.0439) &  25 (0.048)   &    \\
  0.10 &  10   &  77 (0.0497) &  21 (0.030)  &     \\
  0.15 &  10   &  77 (0.0497) &  21 (0.030)  &     \\
  0.20 &  10   &  60 (0.0467) &  18 (0.0380)  &    \\
  0.25 &  10   &   60 (0.0467) &  18( 0.0380)  &     \\
 \hline
    0 &  15    &  152 (0.0402) & 40 (0.029)   & 799 (0.050)\\
 0.05 &  15    & 152 (0.0402) & 40 (0.028) &   \\
 0.10 &  15    & 133 (0.0423) & 36 (0.032)  &    \\
 0.15 &  15    & 116 (0.0494) & 33 (0.048)  &   \\
 0.20 &  15    & 98 (0.0430) &  29 (0.039) &   \\
 0.25 &  15    &  98 (0.0430)  & 29  (0.039)  &   \\
  \hline
    0 &  20    &  211 (0.0413) & 55 (0.040)    & 2145 (0.049) \\
 0.05 &  20    &  192 (0.0463) & 51 (0.036)  &    \\
 0.10 &  20    &  173 (0.0470) & 48 (0.045) &    \\
 0.15 &  20    &  155 (0.0471) & 44 (0.040) &    \\
 0.20 &  20    &  138 (0.0430) & 41  (0.048)   &    \\
 0.25 &  20    &  121 (0.0475) & 37 (0.041) &    \\
    \hline
   0  &  25    &  271 (0.0481)  & 70 (0.032) & 2148 (0.05)\\
 0.05 &  25    &  250 (0.0418)  & 66 (0.037) &     \\
 0.10 &  25    &  231 (0.0452)  & 62 (0.037) &     \\
 0.15 &  25    &  213 (0.0473)  & 59 (0.0465) &     \\
 0.20 &  25    &  177 (0.0462)  & 52 (0.049)  &     \\
 0.25 &  25    &  161 (0.0495)  & 48 (0.041)  &     \\
\hline
\end{tabular}

\end{table}
%
%
Table \ref{cv5table} contains the critical values of the $V_{\rho}$-tests, $M_{\rho}$-tests and Jonckheere-Terpstra test for 5 samples with equal sample sizes $n=5 (5) 25$, and $\rho =0(0.05)0.25$ at 5\% level of significance. The  exact values of the  tail probabilities are equal to the value given  are provided  in parentheses.
\begin{table}
\centering
\caption{Critical values for 5 samples with equal sample sizes $n=5 (5) 25$, and $\rho =0(0.05)0.25$ at near 5\% level of significance}
\label{cv5table}
\begin{tabular}{cc||c|c|c}
$\rho$ &$n$   & $V_{\rho}$&    $M_{\rho}$&    $JT$      \\
\hline
       0 &  5    &  65 (0.0460) &  15 (0.0324) &  160 (0.0492) \\
    0.05 &  5    &  65 (0.0460) &  15 (0.0324) &    \\
    0.10 &  5    &  65 (0.0460) &  15 (0.0324) &    \\
    0.15 &  5    &  65 (0.0460) &  15 (0.0324) &    \\
    0.20 &  5    &  42 (0.0475) &  11 (0.0375)   &    \\
    0.25 &  5    &  42 (0.0475) &  11 (0.0334)  &    \\
\hline
       0 &  10   &  158 (0.0450) &   34  (0.0376)   & 595 (0.0499)\\
    0.05 &  10   &  158 (0.0450) &  34  (0.0376)   &    \\
    0.10 &  10   & 129 (0.0496) &  29  (0.0331)  &     \\
    0.15 &  10    & 129 (0.0496) &  29  (0.0331)  &     \\
    0.20 &  10   & 102 (0.0390) &  25  (0.0334)  &    \\
     0.25 &  10   & 102 (0.0499) &  25  (0.0375) &    \\
 \hline
      0 &  15    & 255 (0.0436)   & 54 (0.0442)   &  1298 (0.0498)\\
   0.05 &  15    & 255 (0.0436)   & 54 (0.0442) &   \\
   0.10 &  15    & 223 (0.0451)   & 48 (0.0352)  &    \\
   0.15 &  15    & 193 (0.0477)   & 44 (0.0459)  &   \\
   0.20 &  15    & 165 (0.0470)   & 40 (0.0479) &   \\
   0.25 &  15    & 165 (0.0489)   & 40 (0.0479)  &   \\
  \hline
   0 &     20    &  353 (0.0428)  & 74 (0.0442)    & 2269 (0.0499) \\
   0.05 &  20    &  320 (0.0464)  & 74 (0.0442)  &    \\
   0.10 &  20    &  288 (0.0439)  & 64 ( 0.0462) &    \\
   0.15 &  20    &  258 (0.0480)  & 59 (0.0459) &    \\
   0.20 &  20    &  229 (0.0472)  & 54 (0.0413)   &    \\
   0.25 &  20    &  203 (0.0499)  & 50 (0.0442) &    \\
    \hline
     0  &  25    & 452 (0.0462)  &  94 (0.0408) &  3500 (0.04989)\\
   0.05 &  25    & 418 (0.0457)  & 89 (0.0489) &     \\
   0.10 &  25    & 385 (0.0439)  & 84 (0.0479) &     \\
   0.15 &  25    & 354 (0.0470)   & 79 (.0485) &     \\
   0.20 &  25    & 295 (0.0479)   & 70 (0.0484)  &     \\
   0.25 &  25    & 267 (0.0483)  &  65 (0.0410)  &     \\
\hline
\end{tabular}

\end{table}

Under $H_0$, the critical region for the $JT$  test for a pre-specified significance level $\alpha$ is  $\left\{JT \ge s_0\right\}$,  where $s_0$ satisfies $ P(JT\ge s_0 | H_0)\le \alpha$. \cite{Jonckheere:54} provides tables of $P(JT \ge s_0|H_0)$ for $k$ samples of equal size $n_{*}$  for: $k$ = 3, $n_{*}\in \left\{2,3,4,5\right\}$; $k$ = 4, $n_{*}\in \left\{2,3,4\right\}$; $k$ = 5, $n_{*}\in \left\{2,3\right\}$ and $k$ = 6, $n_{*}=2$.

\begin{figure}
	\centerline{
		\includegraphics[width= 0.9\textwidth]{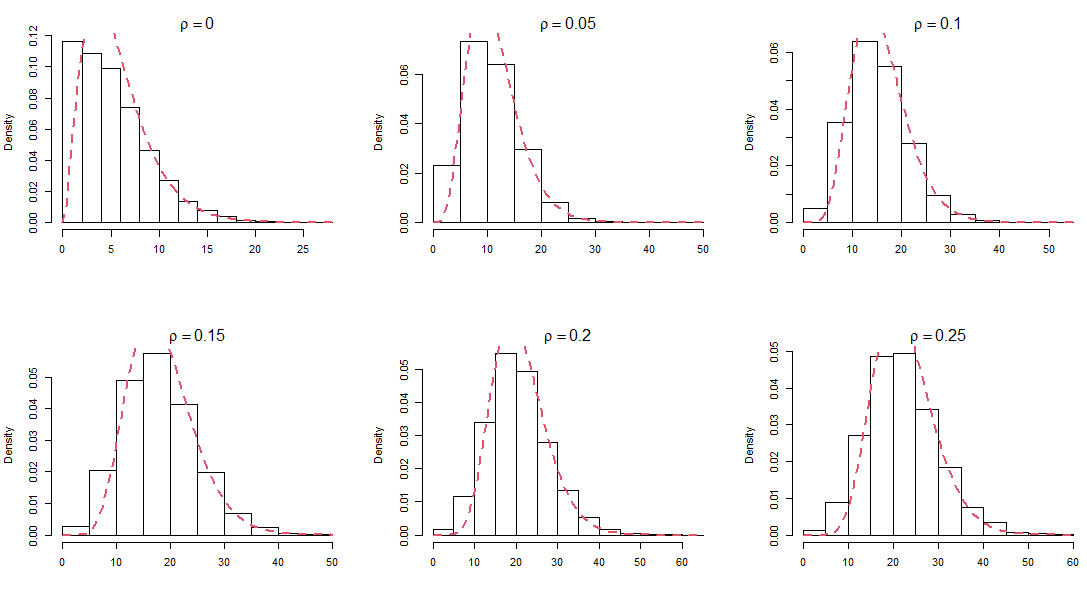}
	}
	\caption{Null distributions of $V_{\rho}$ for 3 samples with equal sample sizes $n_3=20$, and $\rho =0(0.05)0.25$ with $\chi^2$-approximation.}
	\label{null10}
\end{figure}

Since the exact calculations for the null distribution of a rank test statistic are computationally and timely demanding, our empirical experience shows that in some cases the null distributions can be approximated by a continuous $\chi^2$ distribution. Here, we present empirical results for the $k$-sample test statistic $V_{\rho}$ in \eqref{Vrhom}. The results presented below are obtained by 10000 simulation runs.
In Table \ref{gam-appr2}, we provide an example of the exact significance probabilities of  $V_{\rho}$-statistics (at 5\% level), for 3 samples with equal sample size. Figure \ref{null10} shows the distributional shapes  of $V_{\rho}$ for 3 samples with equal sample sizes $n_1=n_2=n_3=20$ and for six choices of the proportion coefficient $\rho$.

Although the $\chi^2$-approximations were designed to have the same expectation as the exact null distributions, it can be noted that their tail behavior also resembles the exact distributions. Hence, for the $V_\rho$ test, we can use the critical regions derived by the $\chi^2$-approximation instead of the exact regions.
\begin{table}
	\centering
	\caption{Comparison of critical values (near 5\% critical values) for 3 samples  with equal sample sizes  $n= 10,15,20$, and $\rho =0(0.05)0.25$ and critical values of $\chi^2$-approximation at 5\% level of significance.}
	\label{gam-appr2}
\begin{tabular}{cc||c|c}
$\rho$ &$n$   & $V_{\rho}$&    $\chi^2$-approx.          \\
\hline
     0 &  10   & 46 (0.0347) & 46.6 (0.05)     \\
  0.05 &  10   & 46 (0.0347) & 46.7  (0.05)      \\
  0.10 &  10   & 38 (0.0490) & 37.3  (0.05)      \\
  0.15 &  10   & 38 (0.0490) & 37.1 (0.05)       \\
  0.20 &  10   & 30 (0.0490) & 29.0 (0.05)        \\
  0.25 &  10   & 30 (0.0490) & 29.4 (0.05)      \\
 \hline
    0 &  15    & 75 (0.0372)   & 76.2 (0.05)     \\
 0.05 &  15    & 75 (0.0372)   & 76.1 (0.05)    \\
 0.10 &  15    &  66 (0.0429)   & 66.32 (0.05)      \\
 0.15 &  15    &  57 (0.0416)   &  57.5 (0.05)     \\
 0.20 &  15    &  49 (0.0469)   & 49.9 (0.05)   \\
 0.25 &  15    &  49 (0.0469)   & 49.4 (0.05)    \\
 \hline
  0 &     20    &  105 (0.0379)  & 105.9   (0.05)     \\
 0.05 &  20    &   95 (0.0370)  & 96.0  (0.05)     \\
 0.10 &  20    &   86 (0.0376)  & 86.8  (0.05)      \\
 0.15 &  20    &   78 (0.0451)  & 78.3  (0.05)      \\
 0.20 &  20    &   69 (0.0387)  & 70.1 (0.05)      \\
 0.25 &  20    &   61 (0.0430)  & 62.2 (0.05)     \\
\hline
\end{tabular}

\end{table}

In the next section, we apply the testing procedures by the proposed $k$-sample rank statistics to the following illustrative example.

\subsection{ Example}\label{3compare}

The data is used to compare the  training and experience of 3 groups of judges, 21 of which were staff members, 23 trainees and 28 were undergraduate psychology majors, see \cite{lehmann1975statistical}.
\begin{table}[h!]
	\centering
\caption{ Accuracies of the judges in terms of percent correctly identified.}
\begin{tabular}{@{}l*{10}{c}}
	Group & \multicolumn{10}{c}{Scores}\\
	\hline
	Staff  & 68.5 & 69.0  & 69.0  & 70.5  & 70.5  & 70.5  & 71.5  & 72.0  & 73.0  & 73.5 \\
	& 73.5  & 74.0  & 74.0  & 74.5  & 75.0  & 75.0  & 75.0  & 75.5  & 76.0  & 76.5  \\
	& 76.5 & & & & & & & & \\
	\hline
	Trainees  & 62.5  & 63.0  & 66.0  & 68.5  & 69.0  & 69.5  & 69.5  & 70.0  & 70.0  & 70.0 \\
	&      70.5  & 70.5  & 71.0  & 71.5  & 71.5  & 71.5  & 73.0  & 73.5  & 74.0  &  74.0 \\
	& 74.0 &  74.5  & 74.5  &   &   &   &   &   &  & \\
	\hline
	Undergraduates  & 58.0  & 60.0  & 64.5  & 65.5  & 66.0  & 66.5  & 68.5  & 68.5  & 69.0  & 69.0 \\
	& 69.0  & 69.0  & 70.0  & 70,5  & 71.0  & 71.5  & 71.5  & 71.5  & 71.5  & 72.0 \\
	& 72.0  & 72.0  & 72.5  & 73.0  & 74.0  & 74.0  & 74.0  & 74.5  &       &
\end{tabular}
\end{table}

\begin{figure}[h!]
\centerline{
 \includegraphics[width=0.5\textwidth]{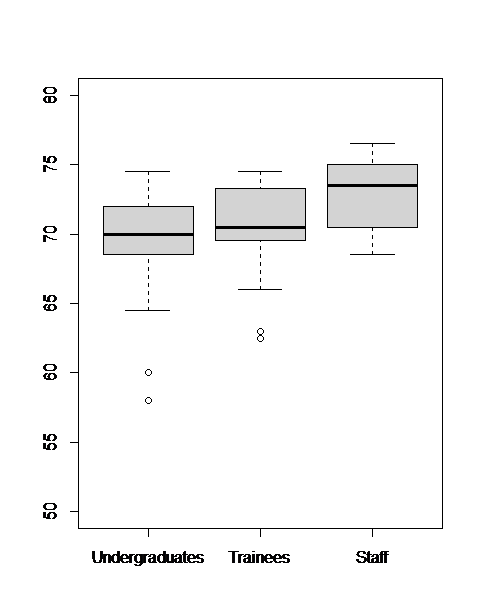}
 }
\caption{Boxplots of the data of Example}%
\label{boxplot}
\end{figure}

In this case, there are given 3 samples with sizes $n_1=21$, $n_2=23$ and $n_3=28$. If training and experience have an effect, the staff members could be expected to be most accurate, the trainees next, and the undergraduates least. Boxplots on Fig. \ref{boxplot} indicate the expected trend ant also indicate that there are possible outliers in the samples.

We compute the statistics $M_{\rho}$ and $V_{\rho}$ for $\rho=0 (0.05) 0.25$.  Table \ref{example1b} contains the values of the statistics, the critical values at 5\% level of significance and the decision about the null hypothesis.

\begin{table}[h!]

\begin{center}
	\caption{ Statistics computation for the Example.}
	\label{example1b}
\begin{tabular}{c|c|c|c||c|c|c}
  \hline
  $\rho$ & $M_{\rho}$ &  c.v. & $H_0$ & $V_{\rho}$ &  c.v. & $H_0$ \\
    \hline
 0    & 40.5 & 44 & reject ($\alpha=0.0243$) & 114.5 & 130 & reject ($\alpha=0.0397$)\\
 0.05 & 36   & 42 & reject ($\alpha=0.0385$) & 106   & 120 & reject ($\alpha=0.0445$)\\
 0.10 & 37   & 39 & reject ($\alpha=0.0354$) & 100.5 & 110 & reject ($\alpha=0.0437$) \\
 0.15 & 25   & 36 & reject ($\alpha=0.0427$) & 80    & 100 & reject ($\alpha=0.0434$)\\
 0.20 & 26   & 33 & reject ($\alpha=0.0332$) & 76.5  & 91  & reject ($\alpha=0.0481$) \\
 0.25 & 27   & 30 & reject ($\alpha= 0.037$) & 69    & 80  & reject ($\alpha=0.0431$)\\
  \hline
\end{tabular}

\end{center}
\end{table}

We have found that all of  of the tests from the two families of statistics reject the null hypothesis of no difference between the three groups at 5\% level of significance.  The observed value of the Jonckheere-Terpstra statistic is $JT = 1164$  with c.v. = 1031 for this data, i.e. JT  also rejects $H_0$ at 5\% level of significance.

\section{Power of the exceedance tests  under Lehmann-type alternatives}\label{power}

The power function of the $M_{\rho}$ test against  the general alternative $H_{A}$ in \eqref{genalt} is given by
\begin{equation}\label{mpow}
   P(M_{\rho}\le s_{\alpha}| H_{A}) = \sum_{j=0}^{s_{\alpha}} P(M_{\rho}=j | H_{A}),
\end{equation}
where $P(M_{\rho}=j | H_{A})$ is the distribution of $M_{\rho}$ under $H_{A}$.

To make meaningful comparison of the power of different $M_{\rho}$ and $V_{\rho}$-tests, we calculated power functions at prescribed exact level  of significance $\alpha$ as follows. First, for any $M_{\rho}$-test, we determine two values $\alpha_L$ and $\alpha_R$ so that
\[
 P(M_{\rho}\le s_{\alpha})=\alpha_L \quad
  \text{and}
  \quad P(M_{\rho}\le s_{\alpha}+1)=\alpha_R,
\]
where $s_{\alpha}$ is given by $P(M_{\rho}\le s_{\alpha}|H_0) \le \alpha$ and therefore, the interval $(\alpha_L,\alpha_R)$ contains the critical level  $\alpha$. Next, we calculate the power values corresponding to the two critical values $s_{\alpha}$ and $s_{\alpha}+1$
\[
  \beta_1 = P(M_{\rho}\le s_{\alpha}|H_{A})
  \quad
  \text{and}
  \quad
  \beta_2 = P(M_{\rho}\le s_{\alpha}+1|H_{A}),
\]
Then, the power of the test at exact level $\alpha$ is estimated by
\[
\beta= \pi\beta_2 + (1-\pi)\beta_1,
\]
where $\pi=\dfrac{\alpha - \alpha_L}{\alpha_R - \alpha_L}$  is the adjusting factor used in the randomized  procedure  in \eqref{rand-prp}.

Similarly, we have adjusted the power values of any $V_{\rho}$-test to the exact level of significance~$\alpha$.

\subsection{Lehmann-type alternative with order statistics }

Now, we will demonstrate the use of   Monte Carlo simulation method for the computation of the power of the exceedance tests $M_{\rho}$ tests and $V_{\rho}$-tests against the Lehmann-type alternatives in \eqref{LE1}.

For this purpose, we have to generate sets of data from the distributions of the order statistics in \eqref{LE1}  and computed the test statistic $M_{\rho}$ and $V_{\rho}$ for each set. Then the power values are estimated by the rejection rates of the null hypothesis for the corresponding $\rho$. Since the distribution of the test statistics $V_{\rho}$ and $M_{\rho}$ do not depend on the underlying distribution $F(x)$ in the alternative hypothesis, we could take the Uniform(0,1) distribution in the calculation.

For a fixed number of samples $k$, we generated 10,000 sets of data from the order statistics $F_1^*(x), \ldots,  F_k^*(x)$  of a Uniform(0,1) distribution and computed the test statistic $V_{\rho}$ and $M_{\rho}$ for each set. The power values were estimated by the rejection rates of the null hypothesis for the corresponding $\rho$.

Comparison of the power of the exceedance tests $M_{\rho}$ and $V_{\rho}$ for 3, 4 and 5 samples  is run using the above Monte Carlo procedure at exact level of significance $\alpha=0.05$.

\begin{figure}
\caption{Power of $V_{\rho}$ and $M_{\rho}$  tests vs. Lehmann alternative for 3 samples}%
\label{ditr0-alt1-3}
\centerline{
 \includegraphics[width=\textwidth]{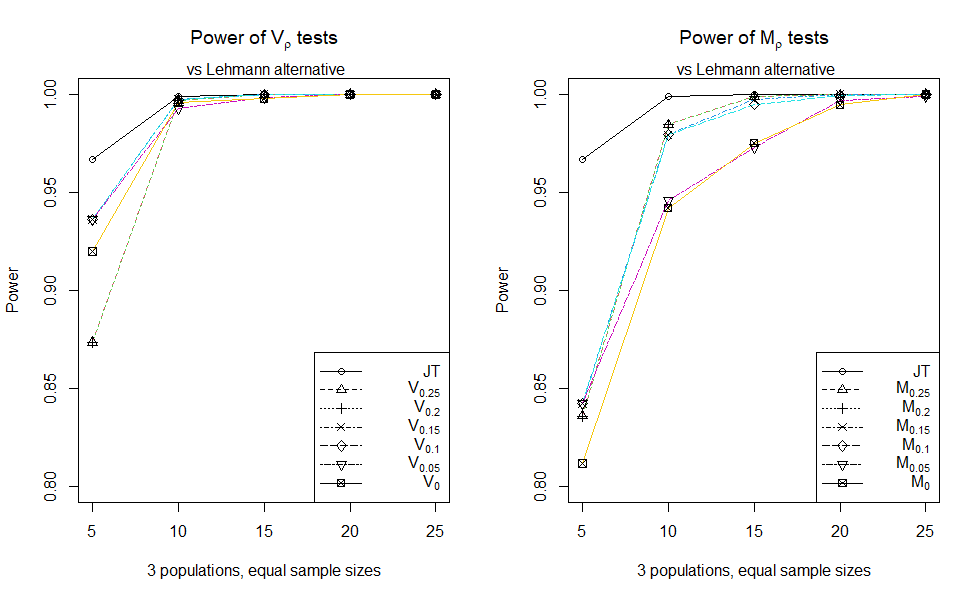}
 }

\end{figure}

\begin{figure}
\caption{Power of $V_{\rho}$ and $M_{\rho}$ tests vs. Lehmann alternative for 4 samples}%
\label{ditr0-alt1-4}
\centerline{
\includegraphics[width=\textwidth]{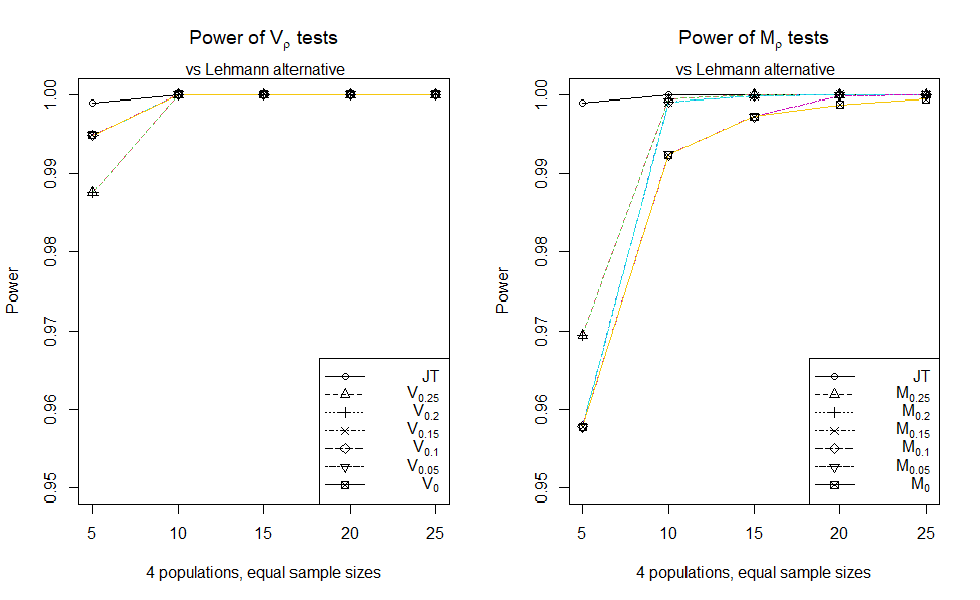}
 }

\end{figure}

\begin{figure}
\caption{Power of $V_{\rho}$ and $M_{\rho}$ tests vs. Lehmann alternativ for 5 samples}%
\label{ditr0-alt1-5}
\centerline{
\includegraphics[width=\textwidth]{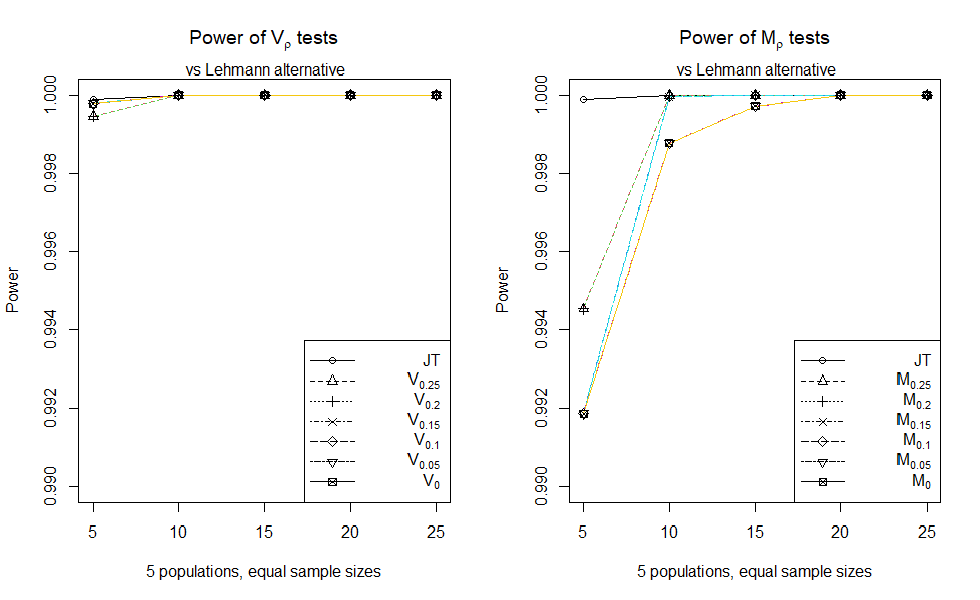}
 }

\end{figure}

Figure~\ref{ditr0-alt1-3} illustrates the  power functions of $V_{\rho}$ and $M_{\rho}$ for 3 samples for $\rho= 0 (0.05) 0.25$  with equal sample sizes $n_i = 5 (5) 25$ ($i=1,2,3$).  The power functions coincide for some tests in some parts since they are equivalent due to the rounding. Figure~\ref{ditr0-alt1-4}  illustrates the  power functions of $V_{\rho}$ and $M_{\rho}$ for 4 samples with $\rho= 0 (0.05) 0.25$ with equal sample sizes $n_i = 5 (5) 25$ ($i=1,2,3,4$).
Figure~\ref{ditr0-alt1-5} illustrates the  power functions of  $V_{\rho}$ and $M_{\rho}$ for 5 samples with $\rho= 0 (0.05) 0.25$ with equal sample sizes $n_i = 5 (5) 25$ ($i=1,2,3,4,5$). The top line in each plot corresponds to the power function of Jonckheere-Terpstra test.

Evidently, for 3 samples the Jonckheere-Terpstra test has the higher power for any sample size between 5 and 25,  while the power of $V_{\rho}$-test  for small $\rho$  ($0\le\rho\le 0.25$)   approaches its level for larger sample sizes. The power of each $M_{\rho}$-test is less than the corresponding  $V_{\rho}$-test but it is also increasing on sample sizes for each $\rho$. We observe similar behaviour of the power functions for 4 and 5 samples. The power of $V_{\rho}$-tests is very high, close to the power of the Jonckheere-Terpstra test. The power of $M_{\rho}$-test is less than the corresponding  $V_{\rho}$-test but still very high.

\subsection{Lehmann-type alternative with power functions}

Next, we compare the power functions against the Lehmann-type alternatives $H_{L_{k,\eta}}$ in \eqref{Leta} with $\eta_i=i$, i.e.
\[
H_{L_{k,i}}: F_i(x) = [F(x)]^{i}, \quad \mbox{where } i= 1, \ldots,k,
\]
where $F(x)$ is an absolutely continuous distribution function.

Comparison of the power of the exceedance tests $M_{\rho}$ and $V_{\rho}$ for $k$ samples  is run using the following Monte Carlo procedure. Since the distribution of the test statistics  $V_{\rho}$ and $M_{\rho}$  do not depend on the underlying distribution $F(x)$ in this alternative hypothesis, we could take the Uniform(0,1) distribution in the calculation. For different number of samples, $k$,  we generate 10,000 sets of data from the distributions  $F_1^*(x)= F(x)$, $F_{2}^*(x) = [F(x)]^2$, \ldots,  $F_k^*(x)=[F(x)]^k$,  and compute the test statistic $M_{\rho}$ and $V_{\rho}$ for each set.

\begin{figure}
\caption{Power of $V_{\rho}$ and $M_{\rho}$ tests for 3 samples vs. Lehmann alternative $F_{n_i}(x) = F(x)^i$}
\label{ditr0-alt2-3}
\centerline{
 \includegraphics[width=\textwidth]{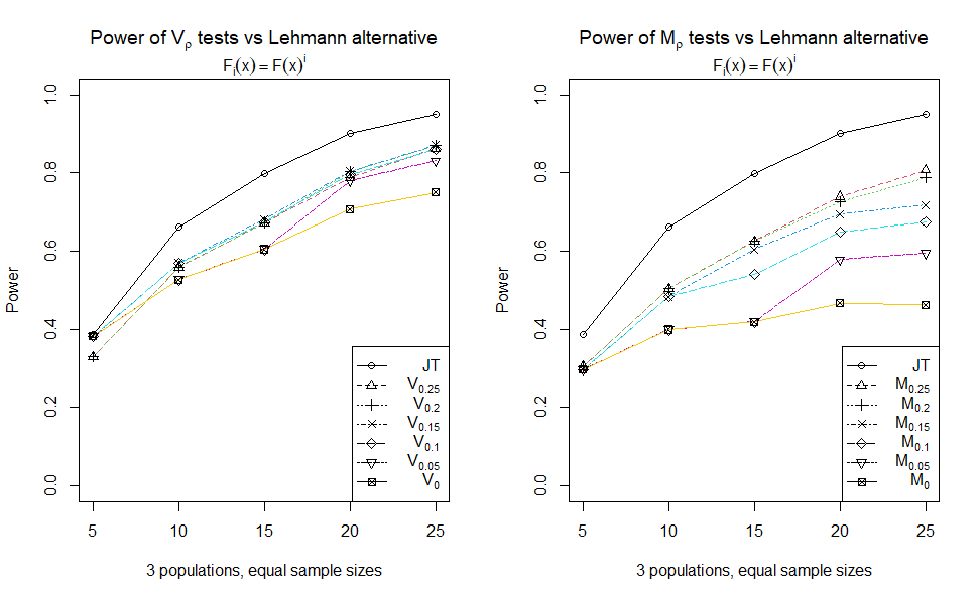}
 }

\end{figure}

\begin{figure}
\caption{Power of $V_{\rho}$ and $M_{\rho}$ tests for 4 samples vs. Lehmann alternative $F_{n_i}(x) = F(x)^i$}
\label{ditr0-alt2-4}
\centerline{
 \includegraphics[width=\textwidth]{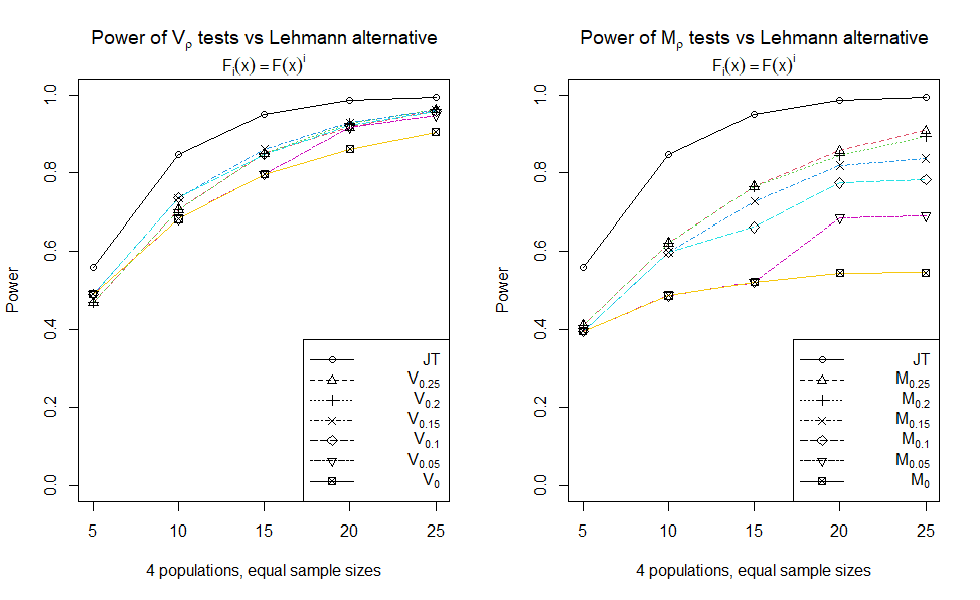}
 }

\end{figure}

\begin{figure}
\caption{Power of $V_{\rho}$ and $M_{\rho}$ tests for 5 samples vs. Lehmann alternative $F_{n_i}(x) = F(x)^i$}
\label{ditr0-alt2-5}
\centerline{
 \includegraphics[width=\textwidth]{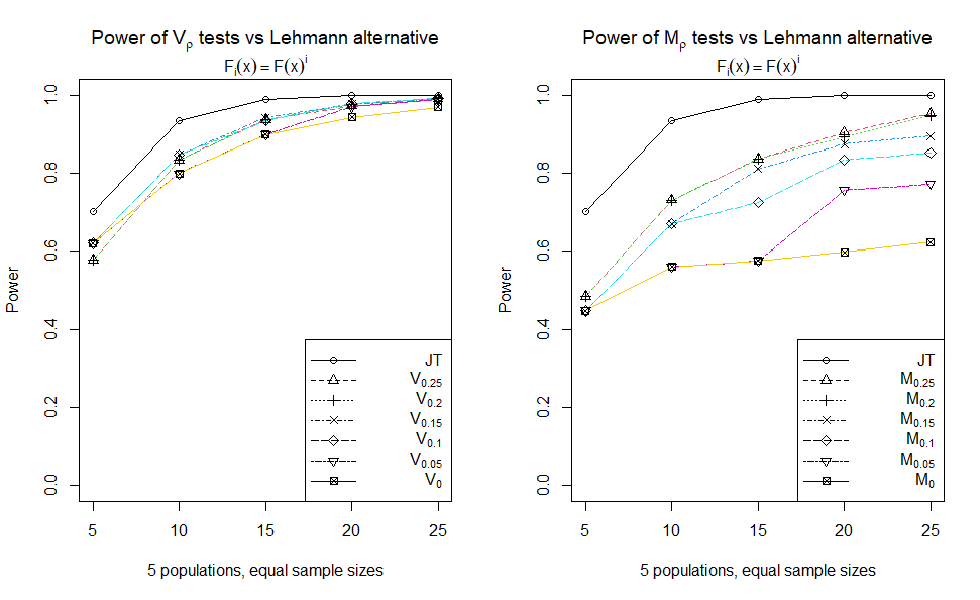}
 }

\end{figure}

Figures \ref{ditr0-alt2-3}-\ref{ditr0-alt2-5} illustrate the  power functions of  $V_{\rho}$ and $M_{\rho}$ with $\rho= 0 (0.05)0.25$ for 3,4 and 5 samples with equal sample sizes $n_i = 5(5)25$. The top line in each plot corresponds to the power function of Jonckheere-Terpstra test.

For this Lehmann-type alternative we could make similar conclusion to the previous case with order statistics. Evidently, for 3 samples the Jonckheere-Terpstra test has the higher power while the power of $V_{\rho}$-tests and $M_{\rho}$-tests for small $\rho$  ($0\le\rho\le 0.25$)   approaches its level for larger sample sizes.

\section{Additional remarks}\label{additionalRemarks}

One of the reasons that Lehmann families of distributions $G(x)=F(x)^\eta$,  for some continuous CDF $F(x)$ and $\eta>0$, proposed by \cite{lehmann53}, are of statistical interest because they lead to an explicit form of the probabilities for the rank orders included in the Neyman-Pearson lemma. Note that the case $\eta\ge1$ is associated with $F(x)\ge G(x)$, while $0<\eta\le1$ corresponds to the alternative $F(x)\le G(x)$. Recently, \cite{balakrishnan2021my} discusses some key distributions and their relation to some of the known models in the distribution theory literature. Substantially, slightly modified Lehmann alternatives in the form $G^{*}(x)= 1 - [1-F(x)]^\eta$ are interesting in life testing problems as they include the class of exponential distributions, the Weibull distributions and distributions with proportional hazard rates, see \cite{lin1992choice} and  \cite{balakrishnan2021my}. Furthermore, the Lehmann family of distributions has the following useful property related to the hazard and reversed hazard rate functions. Suppose the lifetime of a certain item (product) is an absolutely continuous random variable with pdf $f(t)$ and CDF $F(t)$, then the hazard and reversed hazard rate  are defined respectively by
\[
h_F(x) = \frac{f(x)}{1-F(x)}\quad \mbox{and} \quad r_F(x) = \frac{f(x)}{F(x)}; \quad x>0.
\]
If $G(x)$ and $G^{*}(x)$ are another CDFs related to $F(x)$ by $G(x)=[F(x)]^\eta$ and $G^{*}(x)= 1 - [1 - F(x)]^{\eta}$ , respectively, then
\begin{align*}
	r_G(x) &= \frac{g(x)}{G(x)} = \frac{\eta(F(x))^{\eta-1}f(x)}{[F(x)]^\eta} = \eta r_F(x) \quad \mbox{and}\\
	h_G^{*}(x) &= \frac{g^{*}(x)}{1-G^{*}(x)} = \frac{\eta[1 - F(x)]^{\eta-1}f(x)}{[1-F(x)]^\eta} = \eta h_F(x).
\end{align*}
Thus, the parameter $\eta$ has a multiplicative effect on the hazard and reversed hazard rate of $G^{*}(x)$ and $G(x)$, respectively. This feature can be exploit to test the effect on the product lifetime when changing some stress factor in the experiment. In the context of multiple step-stress models, order statistics approach has been discussed by \cite{balakrishnan2012sequential}, while \cite{kateri2015inference, kateri2017hazard} introduced a flexible failure rate-based step-stress accelerated life testing model at various stress levels. In our view, nonparametric inference based on the proposed Lehmann-type alternative with order statistics in \eqref{LE} could be successfully used to compare the observed lifetimes of the experimental units in a sequence of accelerated step-stress levels.

As a final remark, note that in similar settings as in Section \ref{6typesAlternatives}, \cite{heller2006power} studies multi-sample comparison of equality of $k$ samples $H_0 : F_j(x) = F_k(x)$, for $j = 1, \ldots, k-1$, against $H_A:  F_j(x)= \left[F_k(x)\right]^{\eta_j}$, for  $j = 1, \ldots, k-1$, with $\eta_j>1$ for at least one $j$. The distributions in the alternative are not assumed to be ordered. The proposed rank test is related to the Kruskal-Wallis rank test. Heller mention also the extension to the ordered alternative.

\bigskip

\noindent {\bf  Acknowlegdement.} The work of the first author was financed by the European Union - NextGenerationEU, through the National Recovery and Resilience Plan of the Republic of Bulgaria, project No. BG-RRP-2.004-0008. The work of the second author was funded under project KP-06-N52-1 with Bulgarian NSF.

\bibliography{lehmannmulti}

\end{document}